\documentclass[reqno,verbatim]{amsart}
\renewcommand{\Im}{\operatorname{\rm Im}\nolimits}

\newtheorem{thm}{Theorem}[section]
\newtheorem{lemma}{Lemma}[section]

\newtheorem*{lemmab}{Lemma}
\newtheorem{prop}{Proposition}[section]

\theoremstyle{definition}

\theoremstyle{remark}

\def\squarebox#1{\hbox to #1{\hfill\vbox to #1{\vfill}}}

\def \tr {\operatorname{tr}}

\def \supp {\operatorname{supp}}

\def \spec {\operatorname{spec}}

\def \comp {\operatorname{comp}}

\def \S1 {\mbox{S}^{1}}

\def \loc {\operatorname {loc}}

\def \Real {{\mathbb R}}

\def \Complex {\mathbb{C}}

\def \sign {\operatorname{sign}}

\def \Rm {R_-}
\def \Rp {R_+}
\def \Tp {T_+}
\def \Tm {T_-}

\def \mcr {{\mathcal R}}
\def \pro{\Pi}

   \title [Resonances for steplike potentials]
{Resonances for steplike potentials: forward and 
inverse results}
   \author { T. Christiansen}
\thanks{Partially supported by NSF grant 0088922.}

\begin{document}

\begin{abstract}
We consider resonances associated to the operator 
$-\frac{d^2}{dx^2}+V(x)$, where $V(x)=V_+$ if $x>x_M$ and $V(x)=V_-$ if 
$x<-x_M$, with $V_+\not = V_-$.  We obtain asymptotics of the 
resonance-counting function
in several regions.  Moreover, we show that in several situations, the 
resonances, $V_+$, and $V_-$ determine $V$ uniquely up to translation.

\end{abstract}
\maketitle
\section{Introduction}

Suppose $V\in L^{\infty}(\Real;\Real)$ is ``steplike'' in the sense that
there
exist constants $x_M$, $V_+$, $V_-$, $V_+\not = V_-$ so that 
$V(x)=V_{\pm}$ when $\pm x>x_M$. 
We will study the resonances associated with the operator 
$-\frac{d^2}{dx^2}+V$.  The resolvent $R(z)=(-\frac{d^2}{dx^2}+V-z)^{-1}$
is bounded on $L^2(\Real)$ for all but a finite number of $z\in \Complex
\setminus \Real_+$.  Considered as an operator from $L^2_{\comp}(\Real)$
to $ H^2_{\loc}(\Real)$, it has a meromorphic continuation to the 
Riemann surface $\hat{Z}$ which is the minimal Riemann surface on which 
$(z-V_+)^{1/2}$ and $(z-V_-)^{1/2}$ are single-valued holomorphic functions.
This follows from the same techniques as used to study manifolds with 
cylindrical ends
\cite{guce,tapsit}, although this case is much simpler (see Section \ref{s:est}
for an explanation).
  In fact, one motivation
for studying steplike potentials is the similarity their scattering theory
has to the scattering theory for manifolds with cylindrical ends, and a hope
that an understanding of the steplike case will shed some light on what can
happen in the cylindrical ends case.

Choose a $\chi \in C_c^{\infty}(\Real)$ so that $\chi=1$ on the support
of $V'$.  Then the poles of $\chi R(z)\chi$ are independent of the choice
of $\chi$ with this property and we set
$$\mcr = \{z_j\in \hat{Z}:\; \chi R(z)\chi \; \text{has a pole at}\; z_j\}.$$ 
We list the elements with their multiplicities.

Let $\pro: \hat{Z}\rightarrow \Complex$ be the natural projection operator
and let $r_{\pm}(z)=(z-V_{\pm})^{1/2}.$
If $V_+\not = V_-$, $\hat{Z}$ is a four-fold cover of the plane.  We identify
the part of $\hat{Z}$ with $\Im r_+(z)>0$ and
$\Im r_-(z)>0$ as the ``physical sheet,''
that is, the sheet of $\hat{Z}$ on which $R(z)$ is a bounded operator on
$L^2(\Real)$.  

Poles of $R(z)$ are called resonances.  They are in some sense analogous
to eigenvalues and can, in many settings, be viewed as corresponding to 
decaying waves.  See \cite{vodevsurvey, zwsurvey} for an introduction to
resonances and a survey of some results
on their distribution, or \cite{ZwS} for a light-hearted
introduction.

Although resonances have been studied in many settings, there are relatively
few cases where the asymptotics of the resonance counting function are known.
Many of these are in some sense one-dimensional or degenerate, cf. 
\cite{froese, pac,regge,sjrfb, vdegen, zw1}.  
Our first theorem is 
\begin{thm}\label{thm:asym1}
  Suppose $V_+\not = V_-$, and the convex hull of the support of 
$V'$ is $[a,b]$.  Then 
$$\#\{ z_j \in \mcr: |\pro(z_j)|\leq r^2,\; 
\Im r_+(z_j)<0,\; \text{and}\; \Im r_{-}(z_j)<0\}
= \frac{2}{\pi} (b-a)r+o(r).$$
\end{thm}
In case $V_+=V_-$, this theorem is due to \cite{regge,zw1}.  Our proof has 
more in common with the proof of \cite{froese}, though we have to make 
adaptations due to the nature of 
$\hat{Z}$.  In addition, rather than study the zeros of a function
holomorphic on all of $\hat{Z}$, we study the zeros of functions meromorphic 
on $\hat{Z}$, in a region where they are holomorphic except at a finite
number of points.

The next theorem follows by similar techniques, although we must make
stronger assumptions on the perturbation.
\begin{thm}\label{thm:asym2}
Suppose $V(x)=V_+H(x-\beta)+V_-H(\beta-x)+p(x)$, with 
$p\in C^0_{\comp}(\Real)$ and
$p'\in L^1(\Real)$, and the convex hull of the support of $p$ is
$[-b_1,b_1]$.  Then, if $b_1>\beta$,
$$\#\{ z_j\in \mcr:\; |\pro(z_j)|\leq r^2,\; \Im r_+(z_j)<0\; \text{and}\;
\Im r_-(z_j)>0\} = \frac{2(b_1-\beta)}{\pi}r +o(r)$$
and, if $b_1>-\beta$,
$$\#\{ z_j\in \mcr:\; |\pro(z_j)|\leq r^2,\; \Im r_+(z_j)>0\; \text{and}\;
\Im r_-(z_j)<0\} = \frac{2(\beta+b_1)}{\pi}r +o(r).$$
\end{thm}
Here and elsewhere $H$ is the Heaviside function.
If $\beta\in (-b_1,b_1)$, then, we see that the distribution 
of poles on these 
sheets depends rather heavily on the value of $\beta$.  It therefore seems
unlikely that on one of these sheets (say, $\Im r_+>0$,
$\Im r_-<0$) an asymptotic formula holds in general.  However, 
we can prove the weaker
\begin{thm}\label{thm:asym3}  Suppose $V_+\not = V_-$,
and the convex hull of the support of $V'$ is $[a,b]$.
Then
\begin{multline*}
\#\{ z_j\in \mcr:\; |\pro(z_j)|\leq r^2,\; \Im r_+(z_j)\leq 0\; \text{and}\;
\Im r_-(z_j)\geq 0\} \\ + 
\#\{ z_j\in \mcr:\; |\pro(z_j)|\leq r^2,\; \Im r_+(z_j)>0\; \text{and}\;
\Im r_-(z_j)<0\} = \frac{2}{\pi}(b-a)r +o(r).
\end{multline*}
\end{thm}

Zworski \cite{zwinverse} and Korotyaev \cite{k-ir} showed that
 resonances determine 
a compactly supported, even potential in one dimension
if $0$ is not a pole of the resolvent, and that there are two 
such potentials with the same resonances if $0$ is a pole of the 
resolvent.  The papers \cite{b-k-w} and
\cite{k-hl} have studied
a similar
 question for compactly supported potentials on a half-line.  Here we 
give an example of our inverse results for steplike potentials.  Further 
results can be found in Section \ref{s:inverse}.
\begin{prop}\label{prop:inverse1}
Suppose we know all the poles of $R(z)$, $V_+$, $V_-$ (with $V_+\not = V_-$),
and know that our potential $V$ has $V'$ compactly supported.  If, 
in addition, $\pro^{-1}(V_+)\cap \mcr \not = \emptyset$, then
$V$ is uniquely determined up to translation.
\end{prop} 
It is clear that the qualifier ``up to translation'' cannot be removed,
as $V(x)$ and $V(x-c)$, $c\in \Real$, have the same resonances.  The 
techniques of the proof are similar to those  Zworski \cite{zwinverse}
used in recovering an even, compactly supported potential on 
$\Real$.  This problem
is more difficult because on $\hat{Z}$ we do not have the simple 
Weierstrass factorization that we have on $\Complex$.  On the other
hand, because $V_+\not = V_-$,
some ``symmetry'' is broken
and we have more identities at our disposal 
than when $V_+=V_-$.  A related issue is that $\hat{Z}$ is ``bigger''
than $\Complex$, and thus there are ``more'' resonances and they can
carry more information.  

The techniques developed here for the forward problem will be applied in
a future paper to obtain results for the distribution of resonances for
potential scattering on cylinders.

We shall assume throughout that $V_+<V_-$.

\section{Preliminaries from Complex Analysis}
In this section we recall some results and language of complex analysis,
e.g. \cite{levin},
and prove a theorem we shall need on the distribution of zeros of functions
which are ``good'' in a half-plane.

We shall often
work with functions that are holomorphic not in the whole
plane but are holomorphic within an angle $(\theta_1, \theta_2)$.
A function $F$ holomorphic in an angle $(\theta_1,\theta_2)$ is of order
$\rho$ there
 if $$\overline{\lim}_{r\rightarrow \infty}\frac{\ln \ln (\sup_{\theta 
\in \theta_1,\theta_2}|F(re^{-\theta})|)}{r}=\rho.$$
A function of order $\rho$ in the angle $(\theta_1,\theta_2)$ is of type
$\tau$ there if 
$$\overline{\lim} _{r\rightarrow \infty} \frac{\ln 
\sup_{\theta\in (\theta_1,\theta_2)}|F(re^{i\theta})|}{r^{\rho}}=\tau.$$
 A function of order
$1$ and type $\tau<\infty$ (in an angle
$(\theta_1,\theta_2)$) is said to be of exponential type there.  Of course,
$\rho$ and $\tau$ can depend on $(\theta_1, \theta_2)$.

The indicator of a function $F$ holomorphic in an angle $\theta_1 \leq \arg z
\leq \theta_2$
and of order $\rho$ is 
$$h_F(\theta)=\overline{\lim}_{r\rightarrow \infty}\frac{\ln |F(re^{i\theta})|}
{r^{\rho}}.$$
A function $F$ is of completely regular growth within the angle
$(\theta_1,\theta_2)$ if 
$$\lim_{\genfrac{}{}{0pt}{}{r\rightarrow \infty }{ r\not \in E}}\frac{\ln |F(re^{i\theta})|}
{r^{\rho}}=h_F(\theta).
$$
where the set $E\subset \Real_+$ is of zero relative measure and the 
convergence is uniform for $\theta\in (\theta_1,\theta_2)$.

We shall abuse notation slightly and also use the language above for a function
that is holomorphic for $\theta_1 \leq \arg z \leq \theta_2$ and 
$z$ outside of a compact set.

For a function $f$ defined in the upper half plane, define
$$n_{f+}(r)=\#\{ z_j: \; \Im z_j>0\; \text{and}\; f(z_j)=0\}.$$
Here, as elsewhere, all zeros are counted with multiplicities.
\begin{thm}\label{thm:complexanal}
Suppose $f(z)$ is holomorphic in the closed upper half plane $\Im z \geq 0$,
\begin{equation}\label{as1}
|f(z)|\leq Ce^{C|z|}
\end{equation}
 there, $f(0)=1$,
\begin{equation}\label{as3}
\left |\int_{-\infty}^{\infty}\frac{d[\arg f(t)]}{dt}dt\right|<\infty
\end{equation}
and 
\begin{equation}\label{as4}
\left|\int_{-\infty}^{\infty}\frac{\ln |f(t)|}{1+t^2}dt\right|<\infty.
\end{equation}
Then 
$$\lim_{r\rightarrow \infty}\frac{n_{f+}(r)}{r}
=\frac{1}{2\pi} \int_0^{\pi}h_f(\varphi)d\varphi.$$
\end{thm}  A related result is Lemma 3.2 of \cite{froese2}.
\begin{proof}
We first adopt some arguments of \cite[Chapter III, Section 2]{levin};
a similar argument is used to prove Lemma 6.1 of \cite{froese2}.  By
the principle of the argument,
$$2\pi n_{f+}(r)=\int_{-r}^r\frac{d[\arg f(t)]}{dt}dt +\int_0^{\pi}
\frac{d[\arg f(r e^{i\varphi})]}{r d\varphi}r d\varphi.$$
We use the Cauchy-Riemann equations to rewrite this as
$$2\pi n_{f+}(r)=\int_{-r}^r\frac{d[\arg f(t)]}{dt}dt+ r \int_0^{\pi}
\frac{d}{dt}\{\ln |f(te^{i\varphi})|\}_{| t=r}d\varphi.$$
Dividing both sides of the equation by $2\pi r$, and integrating from
$0$ to $r$, results in
$$\int_0^r \frac{n_{f+}(\rho)}{\rho}d\rho
= \frac{1}{2\pi}\int_0^r \frac{1}{\rho}\int_{-\rho}^{\rho}
\frac{d[\arg f(t)]}{dt}dt d\rho +\frac{1}{2\pi}\int_0^{\pi}
\ln |f(re^{i\varphi})|d\varphi.$$
By (\ref{as3}), 
$$\left |\frac{1}{2\pi}\int_0^r \frac{1}{\rho}\int_{-\rho}^{\rho}
\frac{d[\arg f(t)]}{dt}dt d\rho \right| \leq C \ln r,$$
so that we can, given $\epsilon>0$, choose $r_{\epsilon}$ so
that 
$$\left|\frac{1}{r}
\int_0^r \frac{1}{\rho}\int_{-\rho}^{\rho}
\frac{d[\arg f(t)]}{dt}dt d\rho \right| \leq \frac{\epsilon}{6}$$
for all $r>r_{\epsilon}$.

Now we roughly follow the proof of 
\cite[Theorem 3, Chapter III, Section 3]{levin}.  By \cite[Theorem 6, Chapter
V, Section 4]{levin}, (\ref{as1}) and (\ref{as3}) imply that $f(z)$ is a 
function of class A and completely regular growth for $0<\arg z<\pi$,
and $h_f(\theta)=k\sin \theta$.  Because $f$ is of completely regular
growth, if $r$ does not belong to the exceptional set $E$ for the function
$f$, then for a sufficiently large $r_{\epsilon}$ and all $r>r_{\epsilon}$ 
we have
$$
\left|\int_0^{\pi}\frac{\ln |f(re^{i\varphi})|}{r}d\varphi
-\int_0^{\pi}h_f(\varphi)d\varphi\right|<\epsilon/6.$$
Then we have, for $r>r_{\epsilon}$ not in the exceptional set,
$$\left|\frac{1}{r}\int_0^{r}\frac{n_{f+}(\rho)}
{\rho}d\rho -\frac{1}{2\pi}\int_0^{\pi}
h_f(\varphi)d\varphi\right|<\epsilon/3.$$
Then
$$\lim_{r\rightarrow \infty,\; r\not \in E}
\frac{1}{r}\int_0^{r}\frac{n_{f+}(\rho)}{\rho}d\rho=\frac{1}{2\pi}
\int_0^{\pi}
h_f(\varphi)d\varphi.$$
In the same manner as \cite{levin}, proof of Theorem 3, Chapter III, we can 
show that the limit holds without the restriction that $r\not \in E$.

Again following \cite{levin}, we have by the monotonicity of $n_{f+}(r)$ 
that for
$k>1$, 
$$\frac{1}{r}n_{f+}(r)\ln k \leq \frac{1}{r}\int_r^{kr}\frac{n_{f+}(t)}{t}dt
=\frac{1}{r}\int_0^{kr}\frac{n_{f+}(t)}{t}dt-\frac{1}{r}\int_0^{r}
\frac{n_{f+}(t)}{t}dt,$$ giving us
$$\lim \sup_{r \rightarrow \infty}\frac{n_{f+}(r)}{r} 
\leq \frac{k-1}{\ln k}\frac{1}{2\pi}\int_0^{\pi}h_f(\varphi)d\varphi$$
for any $k>1$.  Taking the limit as $k\downarrow 1$ we obtain
$$\limsup_{r\rightarrow \infty}\frac{n_{f+}(r)}{r} \leq \frac{1}{2\pi}\int_0^{\pi}h_f(\varphi)d\varphi.$$
By a similar argument (taking $0<k<1$), we obtain
$$\liminf _{r \rightarrow \infty}\frac{n_{f+}(r)}{r} \geq 
\frac{1}{2\pi}\int_0^{\pi}h_f(\varphi)d\varphi.$$
This proves the theorem.
\end{proof}

We shall use another theorem whose proof is very similar to
the previous one.  If $f(z)$ is 
holomorphic in $\{z\in \Complex:\; |z|\geq R_1\}$, let
$$n_{f,R_1}(r)=\# \{ z_j:\; f(z_j)=0,\; R_1<|z_j|\leq r\}.$$
\begin{thm}\label{thm:complexanal2}  For $|z|\geq R_1$, suppose
$f(z)$ is holomorphic,
$|f(z)|\leq Ce^{C|z|}$, and $f$
is of completely regular growth.  Then
$$\lim_{r\rightarrow \infty}\frac{n_{f,R_1}(r)}{r}=\frac{1}{2\pi}\int_0^{2\pi}
h_f(\varphi)d\varphi.$$
\end{thm}
\begin{proof}  As in the proof of the previous theorem, the principle of the 
argument followed by an application of the Cauchy-Riemann equations 
gives us, 
for $r>R_1$, 
$$2\pi n_{f,R_1}(r)=r\int_0^{2\pi}
\frac{d}{dt}\{\ln |f(te^{i\varphi})|\}_{| t=r}d\varphi
-R_1 \int_0^{2\pi}
\frac{d}{dt}\{\ln |f(te^{i\varphi})|\}_{| t=R_1}d\varphi.$$
Dividing both sides of the equation by $2\pi r$ and integrating from
$R_1$ to $r$, results in 
\begin{multline} \int_{R_1}^r \frac{n_{f,R_1}(\rho)}{\rho}d\rho 
= \frac{1}{2\pi}\int_0^{2\pi}
\ln |f(re^{i\varphi})| d\varphi
-\frac{1}{2\pi}\int_0^{2\pi}
\ln |f(R_1e^{i\varphi})| d\varphi \\ -\frac{1}{2\pi}R_1 \log(r/R_1)
\int_0^{2\pi}\frac{d}{dt}\{\ln |f(te^{i\varphi})|\}_{| t=R_1}d\varphi
.\end{multline}
Then the proof follows the proof of the previous theorem.
\end{proof}

\section{Elementary scattering theory of $-\frac{d^2}{dx^2}+V$}
\label{s:est}
For $z\in \hat{Z}$, let $r_{\pm}(z)=(z-V_{\pm})^{1/2}$.  To a point
$z\in \hat{Z}$ we may associate the two roots $r_{\pm}(z)$, and (if we are
consistent) we can determine a mapping $\hat{Z}\rightarrow \hat{Z}$ by
its action on $r_+(z)$ and $r_-(z)$.  We define
three maps $\omega_+$, $\omega_-$, and
$\omega_{+-}:\hat{Z}\rightarrow \hat{Z}$, by 
\begin{align*} r_{\pm}(\omega_{\pm}(z))& =-r_{\pm}(z)\\
r_{\mp}(\omega_{\pm}(z))& =r_{\mp}(z)
\end{align*}
and 
$$r_{\pm}(\omega_{+-}(z))=-r_{\pm}(z).$$
We shall return to these mappings shortly.

We define two functions $\phi_{\pm}(x,z)$ with the following properties:
$$(-\frac{d^2}{dx^2}+V(x)-z)\phi_{\pm}(x,z)=0$$
and, 
\begin{equation*}
\phi_+(x,z)= \left\{
\begin{array}{ll}
  \Tm(z)e^{-ir_-(z)x}\;  & \text{when }\; x\rightarrow -\infty\\
e^{-ir_+(z)x}+\Rm(z)e^{ir_+(z)x} & \text{when }\; x \rightarrow
\infty 
\end{array}\right.
\end{equation*}
and 
\begin{equation*}
\phi_-(x,z)= \left\{
\begin{array}{ll} e^{ir_-(z)x}+\Rp(z)e^{-ir_-(z)x} & \text{when }\;
 x \rightarrow
-\infty \\
  \Tp(z)e^{ir_+(z)x}\;  & \text{when }\; x\rightarrow \infty.
\end{array}\right.
\end{equation*}

Recall that we assume that $V_+<V_-$.  The following identities
are well known when $\pi(z)
\in \Real$ (e.g. \cite[Theorem 2.1]{c-k}) and can be extended to 
$\hat{Z}$ using the meromorphy of $R_{\pm}$, $T_{\pm}$.
  We have translated them to statements on $\hat{Z}$ using
the mappings $\omega_{\pm}$, $\omega_{+-}$ defined above to obtain
\begin{align}
r_-(z)\Tm(z)&=r_+(z)\Tp(z)\label{eq:T+-}\\
\Rm(\omega_+(z))\Rm(z)&=1\label{eq:r+}\\
T_{\pm}(z)R_{\pm}(\omega_{\mp}(z))& =T_{\pm}(\omega_{\mp}(z))\label{eq:rt+}\\
-r_-(z)\Tm(z)\Rp(\omega_{+-}(z))&=r_+(z)\Rm(z)\Tp(\omega_{+-}(z))
\label{eq:rt}\\
\Tm(\omega_{+-}(z))\Tp(z)+\Rm(\omega_{+-}(z))\Rm(z)&=1 \label{eq:norm}\\
\Tm(\omega_{+-}(z))\Tp(z)+\Rp(\omega_{+-}(z))\Rp(z)&=1\label{eq:norm2}\\
\Rp(\omega_{-}(z))\Rp(z)&=1\label{eq:r-}\\
-r_-(z)\Tm(z)\Tm(\omega_-(z))&=-r_+(z)\Rm(\omega_-(z))+r_+(z)\Rm(z)
\label{eq:rdiff}\\
r_+(z)\Tp(z)\Tp(\omega_+(z))& = r_-(z)\Rp(\omega_+(z))-r_-(z)\Rp(z)
\label{eq:rdiff2}.
\end{align}

We shall also want a family of reference operators, which we define below.  For
$\beta \in \Real$, define the potential $V_{\beta}$ to be
$$V_{\beta}(x)= \left\{ \begin{array}{ll} 
V_- & \text{if}\; x<\beta\\
V_+ & \text{if}\; x>\beta. 
\end{array} 
\right.
$$ 
We define the corresponding set of generalized eigenfunctions,
 $\phi_{\pm, \beta}(x,z)$, which in this case can be given explicitly:
\begin{equation}\label{eq:phi+b}
\phi_{+,\beta}(x,z) = \left\{ \begin{array}{ll}\frac{2r_+(z)}{r_+(z)+r_-(z)}
e^{i\beta(r_-(z)-r_+(z))}e^{-ir_-(z)x} & \text{if}\; x<\beta \\
e^{-ir_+(z)x}+\frac{r_+(z)-r_-(z)}{r_+(z)+r_-(z)}e^{-2ir_+(z)\beta}e^{ir_+(z)x}
& \text{if} \; x>\beta
\end{array}\right.
\end{equation}
and 
\begin{equation}\label{eq:phi-b}
\phi_{-,\beta}(x,z)= \left\{ 
\begin{array}{ll} e^{ir_-(z)x}+\frac{r_-(z)-r_+(z)}{r_-(z)+r_+(z)}e^{2ir_-(z)\beta}e^{-ir_-(z)x} & \text{if}\; x<\beta \\
\frac{2r_-(z)}{r_+(z)+r_-(z)}e^{i\beta(r_-(z)-r_+(z))}e^{ir_+(z)x} & 
\text{if}\; x>\beta.
\end{array}\right.
\end{equation}
We shall denote
\begin{equation}
R_{\beta}(z)=\left( -\frac{d^2}{dx^2}+V_{\beta}-z\right)^{-1}.
\end{equation}
From the explicit representation of $\phi_{\pm, \beta}$
and the corresponding expression for the Schwartz kernel of
$R_{\beta}(z)$ (cf. (\ref{eq:sk})), it is relatively
easy to see that 
$R_{\beta}(z):L^2_{\comp}(\Real)\rightarrow H^2_{\loc}(\Real)$ has a meromorphic continuation to $\hat{Z}$.
Then, since $R(z)=R_{\beta}(z)(I+(V-V_{\beta})R_{\beta}(z))^{-1}$,
$R(z)$ also has a meromorphic continuation to $\hat{Z}$.

\section{Asymptotics for the number of resonances}

In this section we use knowledge about the relationship between poles of
the resolvent and poles of $R_{\pm}$, $T_{\pm}$ to obtain results on
the asymptotics for the number of poles.  Many or
all of these results could
be obtained by more closely following the methods of \cite{froese}.
However, we shall want some of the intermediate results we obtain here
for the inverse results of Section \ref{s:inverse}.

To bound the number of resonances, we shall use the following proposition.
It is essentially a restating of Proposition 3.3 of \cite{cep} for this 
special case.  We remark that while the proofs of \cite{cep} are 
given for the Laplacian, or Laplacian plus compactly supported potential,
on a manifold with cylindrical ends, they work almost without change for
our setting.
\begin{prop}\label{prop:rsm}  Let $z\in \hat{Z}$, and 
suppose $\pro(z)$ is not in the point spectrum of $\frac{-d^2}{dx^2}+V$.
Then, if $\Im r_{\pm}(z)<0$ but $\Im r_{\mp}(z)>0$, then $z$ is a pole of 
$\chi R \chi$ if and only if $z$ is a pole of $R_{\mp}$, and the multiplicities
coincide.  If $\Im r_-(z)<0$ and $\Im r_+(z)<0$, then $z$ is a pole of 
$\chi R \chi$ if and only if it is a pole of $\Rm \Rp-\Tm\Tp$, and the 
multiplicities coincide.
\end{prop}
In fact, in this simple case we can say more.  Looking at the expansions
for $\phi_{+}$ at infinity, if $\Rm$ has a pole of order $k$ at
$z_0$, then so does $\Tm$.  Then using (\ref{eq:T+-}),
unless $\pi(z_0)$ is $V_+$ or $V_-$, $\Tp$ has a pole of 
the same order, and thus so does $\Rp$.  Suppose $\Im r_-(z)<0$,
$\Im r_+(z)<0$, and $\pro (z)\not \in \spec (-\frac{d^2}{dx^2}+V)$.
Since
$$(R_-R_+(z)-T_+T_-(z))^{-1}=R_-R_+(\omega_{+-}(z))-
T_+T_-(\omega_{+-}(z)),$$
and since $T_{\pm}$, $R_{\pm}$ are regular at $\omega_{+-}(z)$, the 
multiplicity of $z$ as a pole of $R_-R_+-T_-T_+$ is equal to its multiplicity
as a pole of $R_+$.

We recall some results and notation of \cite{cep} that will be helpful in 
extending Proposition \ref{prop:rsm}.  Although the following 
lemma may be well known, one reference is Lemma 3.1 of \cite{cep}.
\begin{lemmab}
Suppose $A(z)$ is a $d\times d$-dimensional meromorphic matrix, 
invertible for some value of $z$.  
Then, near $z_0$, it can be put into
the form 
\begin{equation*}
A(z) = E(z)\left(\sum_{j=1}^p(z-z_0)^{-k_j}P_j + \sum _{j=p+1}^{p'}
(z-z_0)^{l_j}P_j+P_0\right)F(z)
\end{equation*}
where $E(z)$, $F(z)$, and their inverses are holomorphic near $z_0$, and
$P_iP_j=\delta_{ij}P_i$, $\tr P_0= d-p'$, $\tr P_i=1$, $i=1,...,p'$.
The $k_j$ and $l_j$ are, up to rearrangement, uniquely determined.
\end{lemmab}
Using the notation of the lemma, the ``maximal multiplicity'' 
$\mu m$ 
of $A$ at $z_0$ is 
$$\mu m_{z_0}(A)=\sum_{j=1}^p k_j$$
(cf. \cite[Definition 3.2]{cep}).

The following proposition improves on Proposition \ref{prop:rsm}.
\begin{prop}\label{prop:rsm2}
Let $z\in \hat{Z}$.  If $\pro (z)\not = V_+,V_-$, then the multiplicity of 
$z$ as a pole of the cut-off resolvent $\chi R\chi$ is equal to the multiplicity
of $z$ as a pole of $\Tm$, $\Tp$, $\Rm$, or $\Rp$.  If $\pro (z)=V_+$, then
$z$ is a pole of $\chi R\chi$ of order at most one and is a pole of 
$\chi R\chi$ if and only if $\Tm(z)\not = 0$.  If $\pro (z)\in (V_+,V_-]$,
then $z$ is not a pole of $\chi R \chi $ or of $\Tm$, $\Tp$, $\Rm$ or $\Rp$.\
\end{prop}
\begin{proof}  By Proposition
\ref{prop:rsm} and the subsequent remarks, we need only 
consider what happens when $\pro (z)$ lies in the spectrum of 
$-\frac{d^2}{dx^2} +V$.
We will use
\begin{equation}\label{eq:phi}
\phi_{\pm}(x,z(k))= \phi_{\pm,\beta}(x,z(k))-\left(R(z(k))
(V-V_{\beta})\phi_{\pm,\beta}(\cdot,z(k))\right)(x).
\end{equation}  Thus $R_{\pm}$, $T_{\pm}$ cannot have a pole of multiplicity
greater than the multiplicity of the pole of $R$.

First, suppose $z$ corresponds to an eigenvalue; that is, $z$ lies in the 
physical sheet and $\pro (z)$ lies in the point spectrum of $-\frac{d^2}{dx^2}+V$.  Then  
$\chi R \chi$, $T_{\pm}$, $R_{\pm}$ all have poles of order one at $z$.  
Then $T_{\pm}$, $R_{\pm}$ are regular at $\omega_+(z)$, $\omega_-(z)$,
and $\omega_{+-}(z)$.  By \cite[Theorem 3.1]{cep}, then, $\chi R \chi$
is regular at 
$\omega_+(z)$, $\omega_-(z)$,
and $\omega_{+-}(z).$

By \cite[Proposition 6.7]{tapsit} and the fact that there are no 
eigenvalues embedded in the continuous spectrum, if $z\in \mcr$ 
lies on the boundary
of the physical sheet and $\pro (z) \in [V_+,\infty)$, then
  $\pro(z)=V_+$ or $\pro(z)=V_-$.  Moreover, again by Proposition 6.27 of
\cite{tapsit}, if $\pro(z)=V_-$, $z$ is not a pole of the resolvent,
so the entries in the scattering matrix must be regular there.

Because $\Rm(z)\Rm(\omega_+(z))=1$ and $\Rm(\omega_+(z))=\overline 
{R}_-(z)$ when $\pi(z)\in [V_+,V_-]$, $\Rm(z)$, and hence $T_{\pm}(z)$,
$\Rp(z)$ have no poles when $\pro(z)\in (V_+,V_-).$  Then, by \cite[Proposition
6.7]{tapsit} and 
\cite[Theorem 3.1]{cep}, $\chi R \chi$ has no poles in this region.

Next suppose that $\pro(z)\in (V_-,\infty)$ but that $z$ does not lie 
on the boundary of the physical sheet.  Let 
$$S_2(z)=\left( \begin{array}{cc}
R_-(z) & T_+(z) \\
T_-(z) & R_+(z) 
\end{array}
\right). 
$$  By Theorem 3.1 of \cite{cep}, the multiplicity of $z$ as a pole of
$R$ is equal to its ``maximal multiplicity'' as a pole of $S_2$.
Since, for the $z$ we are considering, $(S_2(z))^{-1}=S_2(\omega_{+-}(z))
=\overline{S}_2(z)$, the determinant of $S_2$ is regular for these values 
of $z$.  Therefore, since $S_2$ is a $2\times 2$ matrix, the maximal 
multiplicity of $z$ as a pole of $S_2$ is equal to its multiplicity as
a pole of $R_+$ (or $T_{\pm}$ or $R_-$).

Finally, we consider what happens for $z_0$ such that $\pro(z_0)=V_+$.  If
$z_0$ lies on the boundary of the physical sheet, by \cite[Proposition 6.7]
{tapsit}, $z_0$ is a pole of order
one of $\chi R \chi$ if and only if $\Tm(z_0) \not = 0$.  Notice in this case
$\Tp$ and $\Rp$ have a pole at $z_0$ but are regular at $\omega_-(z_0)$.   
If $\Tp$ is regular at $z_0$, $\chi R \chi$
is regular at $z_0$.  

Now suppose $z_0$ does not lie on the boundary of
the physical sheet and $\pro (z_0)=V_+$.  Then we will use 
$$R(z)-R(\omega_{+-}(z))= \frac{i}{2r_+(z)}\phi_+(z)\otimes \phi_+(\omega_+(z))
+ \frac{i}{2r_-(z)}\phi_-(z) \otimes \phi_-(\omega_{+-}(z)).$$
Here we use the notation
$(f\otimes g)h(x)=f(x)\int g(x')h(x')dx'$.  Note that $\phi_+$ is
holomorphic at $z_0$ and $\omega_{+-}(z_0).$  Considering  
the expansion at infinity of the residue
of the kernel of $R(z)$, we see that $R$ has a pole at $z_0$ only
if $\phi_-$, and thus $\Tp$, has a pole at $z_0$.   
  Thus,
using (\ref{eq:phi}), $z_0 \in \mcr$ if and only if $\Tp$ has a pole
of order $1$ at $z_0$.
\end{proof}

We shall now study the behaviour of $T_{\pm}(z)$, $R_{\pm}(z)$, when
$z$ lies on the closure of 
 the physical sheet ($\Im r_+(z)>0$, $\Im r_-(z)>0$).  
It will be more convenient to introduce 
the variable $k=r_+(z)$ on its closure.  We write $z=z(k)=(k^2+V_+)^{1/2}$
in this region.  When $k_0\in \Real$, we understand $z(k_0)$ to be the element
of $\hat{Z}$ obtained by taking the limit of $z(k)$ as $k\rightarrow k_0$, $\Im k >0$.

The following proposition is related to Lemma 3.3 and Corollary 3.4 of
\cite{froese}.  
\begin{prop}\label{prop:t}
Suppose that $V'$ is compactly supported.  Then,
on the closure of the set with 
$\Im k>0$ and $\Im r_-(z(k))>0,$ $|T_{\pm}(z(k))-1|\leq C/|k|$.  
\end{prop}
\begin{proof}
We shall assume that $\supp V'\subset[-b_1,b_1],$ and that $b_1>0$.

If $\chi \in C_c^{\infty}(\Real)$ is $1$ on the support of 
$V-V_{\beta}$, 
\begin{equation}\label{eq:series}
R(z)\chi
= R_{\beta}(z)\chi
\sum_{0}^{\infty}(-1)^j\left((V-V_{\beta})R_{\beta}(z)\chi\right)^j
\end{equation}
when $z$ is on the closure of
 the physical sheet and $|z|$ is sufficiently large.  
We shall apply this identity and (\ref{eq:phi}) with $\beta=0$.  The
operator
$e^{\pm ikx}R_{\beta}(z(k))e^{\mp ikx}$ has Schwartz kernel given by 
\begin{equation}\label{eq:sk}
\frac{e^{\pm ik(x-x')}}{W_{\beta}(z(k))}
\left( H(x-x')\phi_{-,\beta}(x,z(k))\phi_{+,\beta}(x',z(k))+ H(x'-x)
\phi_{-,\beta}(x',z(k))\phi_{+,\beta}(x,z(k))\right)
\end{equation}
where
\begin{equation}\label{eq:W}
W_{\beta}(z(k))=
W[\phi_{-,\beta}(x,z(k)),\phi_{+,\beta}(x,z(k))]= 
\frac{-4ikr_-(z(k))}{k+r_-(z(k))}e^{i\beta(r_-(z(k))-k)}
\end{equation}
is the Wronskian.

Let $\tilde{\chi}\in C_c^{\infty}(\Real)$
and note that when $z$ is in the closure of the physical sheet,
$r_+(z)=r_-(z)+O(1/|z|^{1/2})$.  Examining (\ref{eq:phi+b}),
(\ref{eq:phi-b}), (\ref{eq:sk}), and (\ref{eq:W}), we see that 
$$\| \tilde{\chi}(x)e^{\pm ik(x-x')}R_0(x,x',z(k))
\tilde{\chi}(x')\|_{L^{\infty}(\Real_x\times \Real_{x'})}
\leq \frac{C}{|k|}$$
 when $k$ is in the closed upper half plane and 
$z$ is in the closure of the physical sheet.  (The constant
$C$ depends on $\tilde{\chi}$, of course.)  Therefore, in this same
region
$$ \| \tilde{\chi} e^{\pm ikx} R_0(z(k))\tilde{\chi}e^{\mp ikx}\|_{L^2(\Real)
\rightarrow L^2(\Real)}\leq \frac{C}{|k|}$$
and, when $|k|$ is sufficiently large,
\begin{equation}\label{eq:sumbound}
\|\sum_{1}^{\infty} (-1)^j\left(e^{\pm ikx}(V-V_0)R_0(z(k)\chi)e^{\mp ikx}
\right)^j \|\leq 
\frac{C}{|k|}.
\end{equation}
Let 
\begin{align}\label{eq:g}
g(x,k)& = \sum_0^{\infty}(-1)^j \left((V-V_0)R_0(z(k))\chi\right)^j(V-V_0)\phi_{+,0}
\\
\nonumber
 & = e^{-ikx}\sum_0^{\infty}(-1)^j(e^{ikx}(V-V_0)R_0(z)\chi e^{-ikx})^j
(V-V_0)(e^{ikx}\phi_{+,0}).
\end{align}
Then $g$ is supported in $[-b_1,b_1]$, and using (\ref{eq:phi+b}) and 
(\ref{eq:sumbound}) we have
\begin{equation}\label{eq:gbound}
\| e^{ik\cdot}g(\cdot,k)\| \leq C.
\end{equation}
Using (\ref{eq:phi}), (\ref{eq:series}), and the explicit expression for
the Schwartz kernel of $R_{\beta}(z)$, we obtain
\begin{align*}
& \Tm(z(k))-\frac{2k}{k+r_-(z(k))}\\ & = 
\frac{-1}{W_0(z(k))}\frac{2k}{k+r_-(z(k))}
\int _{-b_1}^0 \left(e^{ir_-(z(k))x'}+\frac{r_-(z(k))-k}{r_-(z(k))+k}
e^{-ir_-(z(k))}\right)g(x',k) dx'  \\ & \hspace{4mm}\mbox{}
-\frac{1}{W_0(z(k))}\frac{2k}{k+r_-(z(k))}
\int_0^{b_1}\frac{2r_-(z(k))}{k+r_-(z(k))}e^{ir_+x'}g(x',k)dx'.
\end{align*}
Thus $|\Tm(z(k))-1|\leq \frac{C}{|k|}$ for $k$ in the closed upper half
plane, and $z$ in the closure of the physical space.  

The result for $\Tp$ 
follows from the relation (\ref{eq:T+-}).
\end{proof}

We recall a slight modification of Lemma 4.1 of \cite{froese}.
\begin{lemma}\label{l:froese}
Suppose $h\in L^{\infty}(\Real)$ has support contained in $[-1,1]$, but
in no smaller interval.  Suppose $f(x,k)$ is analytic for $k$ in the 
closed upper
half plane, and for real $k$ we have $f(x,k)\in L^2([-1,1]dx, \Real dk)$.
Then $\int e^{\pm ikx}h(x)(1-f(x,k))dx$ has exponential type at least $1$ for
$k$ in the upper half plane.
\end{lemma} 

We will use this lemma to prove
\begin{prop}\label{prop:r}
Suppose that $[-b_1,b_1]$ is the convex hull of the support of $V'$.  Then
$R_{\pm}(z(k))$ are, for $\Im k>0$,
functions of exponential type and completely regular growth.  Moreover,
$R_{\pm}(z(k))$ is of type $2b_1$ in this region.  For 
$t\in \Real$ and fixed $\alpha \geq 0$,
and $z(t+i\alpha)$ in the closure of the physical sheet, $R_{\pm}
(z(t+i\alpha))
= O(|t|^{-1})$. 
\end{prop}
\begin{proof}
The proof of this Proposition resembles that of the previous one, though 
we must use Lemma \ref{l:froese}.

We give the proof for $\Rm$, as the proof for $\Rp$ is similar.
Using the functions $g$ and $W_0$ defined in (\ref{eq:g}) and (\ref{eq:W}),
respectively, we have
$$\Rm(z(k))- \frac{k-r_-(z(k))}{k+r_-(z(k))}= 
\frac{-1}{W_0(z(k))}\frac{2r_-(z(k))}{k+r_-(z(k))}(I_1+I_2+I_3)$$
where 
\begin{align}
I_1(k)& =\int_{-b_1}^{0}\frac{2k}{k+r_-(z(k))}e^{-ir_-(z(k))x'}g(x',k)dx'\\
I_2(k)& = \int_0^{b_1}e^{-ikx'}g(x',k)dx'\\
I_3(k)& = \int_0^{b_1}\frac{k-r_-(z(k))}{k+r_-(z(k))}e^{ikx'}g(x',k)dx'.
\end{align}
Using (\ref{eq:gbound}), $|I_1(k)|\leq C$, $|I_3(k)|\leq C$
when $\Im k \geq 0$ and 
$|k|$ is sufficiently large.  We rewrite
\begin{align}\label{eq:I2}
I_2(k)& = \int_0^{b_1}e^{-2ikx'}(V-V_0)(x')(1+f_1(x',k))dx' + 
\int_0^{b_1}e^{-ikx'}\frac{k-r_-(z(k)}{k+r_-(z(k))}e^{ikx'}(V-V_0)(x')dx'
\end{align}
where
$$
f_1(x',k)=\sum_{1}^{\infty}(-1)^j( e^{ikx'}\chi
R_0(z(k))(V-V_0)e^{-ikx'})^je^{ikx'}\phi_{+,0}.$$
The
 second integral in (\ref{eq:I2}) is clearly bounded.  
Using (\ref{eq:sumbound}), 
the first integral is, for 
large $|k|$, bounded by $Ce^{2b_1\Im k}$.  This shows that
$|\Rm(z(t+i\alpha))|=O(|t|^{-1}).$

The function $f_1(x,k)$ may have poles a finite number of points with
$\Im k\geq 0$ (and $\Im r_-(z(k))>0$).  These poles correspond to eigenvalues
of $-\frac{d^2}{dx^2}+V$.
Using (\ref{eq:sumbound}), for $t$ real, and for all but 
a finite number of $\alpha \geq 0$,
 $f_1(x',t+i\alpha)\in L^2([-b_1,b_1]dx', \Real dt)$.  Now, by a
 shift of variable (using $k-i\alpha$, for $\alpha$ chosen
greater than the largest imaginary part of a pole of $R(z(k))$)
followed by a rescaling, we may apply 
Lemma \ref{l:froese} to see that $I_2(k)$, and thus $\Rm(z(k))$, is of type
$2b_1$.
\end{proof}

\begin{proof}[Proof of Theorem \ref{thm:asym1}]
We continue to use the variable $k=r_+(z)$, working in the closed upper
half plane (with $\Im r_-(z(k))\geq 0$).

Suppose $k_1,k_2,...,k_n$ are the poles, listed with 
multiplicity, of $\phi(k)=(\Tm\Tp-\Rm \Rp)(z(k))$ for $k$ in the closed 
upper half plane.  Let $p(k)=\prod_{j=1}^n(1-k/k_j)$.  Then
$\phi_1(k)=p(k)\phi(k)/p(-k)$ is a holomorphic function in $\Im k>0$, with
a continuous extension to the closed half plane.   Now let 
$\phi_2(k)=\phi_1(k)/\phi_1(0)$ if $\phi_1(0)\not = 0$ and $\phi_2(k)=\phi_1(k)m!/\phi_1^{(m)}(0)
k^m$ if $\phi^{(j)}_1(0)=0$ for $0\leq j<m$ and $\phi_1^{(m)}(0)\not = 0$.


By Proposition \ref{prop:t},
$$(\Tm\Tp -\Rm \Rp)(z(k))=1-\Rm \Rp(z(k))+O(|k|^{-1}).$$
Thus
$$h_{\phi_2(k)}(\varphi)=h_{\Rm R_+(z(k)}(\varphi)=4b_1\sin \varphi$$
using Proposition \ref{prop:r}
and the fact that $\Rm(z(k)),\; \Rp(z(k))$ are of completely regular 
growth in the upper half plane.  Since $R_{\pm}(z(t))=O(|t|^{-1})$
for $t$ in $\Real$, we have
$
|\int_{-\infty}^{\infty}\frac{d[\arg \phi_2(t)]}{dt}dt|<\infty$.  Although
$\phi_2(k)$ is not holomorphic in the {\em closed} upper half plane, it 
is holomorphic in the closed upper half plane except at the points
$k=\pm (V_--V_+)^{1/2}.$  It is continuous at these points (in fact, 
a stronger statement can be made), and that is enough to guarantee that we
may apply Theorem \ref{thm:complexanal} to obtain Theorem \ref{thm:asym1}.
\end{proof}  

We next turn to proving Theorem \ref{thm:asym2}.  For this, we shall
need to better understand the behaviour of $R_{\pm}(z(k))$ when $k\in \Real$.
Recall that $\Im (r_-(z(k)))>0$ when $\Im k>0$, and thus 
the restriction on the sign of $r_-(k)$ ensures that we are on the boundary
of the physical sheet.
\begin{lemma}\label{l:Vnice}
Suppose $V(x)=V_{\beta}+p(x)$, with 
$p\in C^0_{\comp}(\Real)$, and $p'\in L^1(\Real)$,
 and the convex hull of the support of $p$ is
$[-b_1,b_1]$.  Then, for $k\in \Real$,
and $\sign r_-(z(k))=\sign k$, $\Rm(z(k))=
\frac{k-r_-(z(k))}{k+r_-(z(k))}e^{-2ik\beta}+o(1/k^2)$ and $\Rp(z(k))=
\frac{r_-(z(k))-k}{r_-(z(k))+k}e^{2ir_-(z(k))\beta}+o(1/k^2).$
\end{lemma}
\begin{proof}
We give the proof for $\Rm$.  Using (\ref{eq:phi}) and (\ref{eq:series}),
for large $|k|$, $k\in \Real$ we can write
$$\phi_+(x,z(k))-\phi_{+,\beta}(x,z(k))=-
R_{\beta}(z(k))(g_1-g_2+g_3)(x,k)$$
where
\begin{align*}
g_1(x,k)& =(p\phi_{+,\beta})(x,(z(k))\\
g_2(x,k)& = (pR_{\beta}p\phi_{+,\beta})(x,z(k))\\ 
g_3(x,k)& = \sum_{j=2}^{\infty}(-1)^j
((pR_{\beta})^jp\phi_{+,\beta})(x,z(k)).
\end{align*}
Therefore, 
$$\Rm(z(k))-\frac{k-r_-(z(k))}{k+r_-(z(k))}e^{-2ik\beta} 
= \frac{-1}{W_{\beta}(z(k))} \frac{2r_-(z(k))}{k+r_-(z(k))}
e^{i\beta(r_-(z(k))-k)}\int \phi_{+,\beta}(x',z(k))(g_1-g_2+g_3)(x')dx'.$$
Since $\|g_3(\cdot, z(k))\|=O(|k|^{-2})$ for $k\in \Real$, we need only
concern ourselves with the contributions of $g_1$ and $g_2$.

We have
\begin{align*}
& \int \phi_{+,\beta}(x',z(k))g_1(x',k)dx' \\ & =
H(\beta+b_1)
\int_{-b_1}^{\beta} p(x')\frac{4k^2}{(k+r_-(z(k)))^2}e^{2i\beta(r_-(z(k))-k)}
e^{-2ir_-(z(k))x'}dx' \\ & \hspace{3mm}
\mbox{}+H(b_1-\beta)\int_{\beta}^{b_1}p(x')\left(
e^{-2ikx'}+\frac{2(k-r_-(z(k))}{k+r_-(z(k))}e^{-2ik\beta}
+\left(\frac{k-r_-(z(k))}{k+r_-(z(k))}\right)^2
e^{-4i\beta k}e^{2ikx'}\right)dx'.
\end{align*}
We will use the fact that $r_-(z(k))=k+O(|k|^{-1})$.  Using integration
 by parts
in the first integral and in the first term of the second one, we obtain
\begin{multline}
p(\beta)\left( \frac{4k^2H(\beta+b_1)}{-2ir_-(z(k))(k+r_-(z(k))^2}
+\frac{H(b_1-\beta)}{2ik}\right)e^{-2i\beta k}
\\
+ \frac{H(\beta+b_1)}
{2ir_-(z(k))}\int_{-b_1}^{\beta}\frac{p'(x')4k^2}{(k+r_-(z(k)))^2}
e^{2i\beta(r_-(z(k)-k)x'}e^{-2ir_-(z(k))x'}
dx' \\ \mbox{}+\frac{H(b_1-\beta)}{2ik}
\int_{\beta}^{b_1}p'(x')e^{-2ikx'}dx' +O(|k|^{-2}). \label{eq:t1}
\end{multline}
The first term is $O(|k|^{-2})$ using the fact that $r_-(z(k))=k+O(|k|^{-1})$.
If $h\in L^1(\Real)$, then the Fourier transform $\hat{h}$ satisfies
$\hat{h}(k)=o(1)$.   Therefore, the two integrals in 
(\ref{eq:t1}) contribute terms which are $o(|k|^{-1})$, so that
$\int \phi_{+,\beta}(x',z(k))g_1(x',k)
=o(|k|^{-1})$.

Next we consider
$\int \phi_{+,\beta}(x',z(k))g_2(x',k)dx'.$  Note that 
$\|g_2(\cdot,z(k))\| =O(|k|^{-1})$ and $g_2$ has compact support, so that 
$$\int (\phi_{+,\beta}-e^{-ikx'})(x',z(k))g_2(x',k)dx'=O(|k|^{-2}).$$
Similarly, 
$\|g_2(x,k)- (pR_{\beta}pe^{-ik\cdot})(x,z(k))\|_{L^2(\Real_x)}
=O(|k|^{-2})$ so that
we may further simplify our calculations.  Moreover, using the 
explicit formula for the resolvent in terms of $\phi_{\pm,\beta}$,
$$
\left \| (pR_{\beta}p(\cdot)e^{-ik\cdot})(x,z(k)) 
- \frac{1}{W_{\beta}(z(k))}p(x)\int  e^{ik|x-x'|}p(x')e^{-ikx'}dx' 
\right\|_
{L^2(\Real_x)}=O(|k|^{-2}).$$
Thus
\begin{equation}\label{eq:g2}
\int \phi_{+,\beta}(x',z(k))g_2(x',k)dx'
= \frac{1}{W_{\beta}(z(k))}\int\int p(x')e^{-i(x'+x'')k}e^{ik|x''-x'|}p(x'')dx'dx''
+O(|k|^{-2}).
\end{equation}
To show that (\ref{eq:g2}) is $o(|k|^{-1})$, we can integrate by parts in
$x''$ when $x''<x'$ and in $x'$ when $x''>x'$, thus finishing the proof of 
the lemma.
\end{proof}

We shall use the previous lemma, Theorem \ref{thm:complexanal},
Proposition \ref{prop:rsm}, and 
Proposition \ref{prop:r} to prove Theorem \ref{thm:asym2}.
\begin{proof}[Proof of Theorem \ref{thm:asym2}]  If $\Rm(k)$ has no
poles in $\Im k>0$ 
we shall apply Theorem \ref{thm:complexanal} to the function
$F(k)=\Rm(z(k))e^{2ik\beta}/\Rm(0)$. 
 Note that for $k\in \Real$, $F(k)=(V_--V_+)(4k^2R_-(0))^{-1}+o(|k|^{-2})$
as $|k|\rightarrow \infty.$  Moreover, by 
Proposition \ref{prop:r}, $h_{F}(\varphi)=2(b_1-\beta)\sin
\varphi$.  Then we may apply Theorem 
\ref{thm:complexanal} to find that the number of zeros of $F(k)$ in the
upper half plane with norm less than $r$ is given by 
$2(b_1-\beta)(\pi)^{-1}r+o(r).$  Then the first part of the theorem 
follows by applying Proposition \ref{prop:rsm}.

If $\Rm$ has any poles in the upper half-plane, it has only finitely many,
and these may be handled as in the proof of Theorem \ref{thm:asym1}.
The second part of the theorem follows in an analogous way, using the
estimates on $\Rp$.
\end{proof}

We now give the proof of the last of our principle forward results, 
Theorem \ref{thm:asym3}.
\begin{proof}[Proof of Theorem \ref{thm:asym3}]

By translation, we may assume that the convex hull of the support of $V'$
is $[-b_1,b_1]$.

As before, use the coordinate $k=r_+(z)$ to describe points in the closure
of the physical sheet of $\hat{Z}$.  Consider the function 
$F(k)$ defined for $k\in \Complex \setminus [-\sqrt{V_--V_+}, \sqrt{V_--V_+}]$
by
$$F(k)=\left\{ \begin{array}{ll}
\frac{\Rm(z(k))}{\Tm(z(k))} & \text{if}\; \Im k \geq 0,\;
 k \not \in [-\sqrt{V_--V_+}, \sqrt{V_--V_+}]\\ 
-\frac{\Rp(z(-k))}{\Tm(z(-k))}& \text{if}\; \Im k<0.
\end{array}
\right.
$$  Note that $F(k)$ depends on the values of $R_{\pm}$, $\Tm$ only
on the closure of the
physical sheet, although $F$ is defined on $\Complex 
\setminus [-\sqrt{V_--V_+}, \sqrt{V_--V_+}]$.
We remark that $\Tm(z(k))$ is nonzero in this region, and that if 
$T_-$ has a pole, then $R_{\pm}$ has a pole of the same order at the same 
place.
The function
$\Rm(z(k))/\Tm(z(k))$ is holomorphic on $\{k: \Im k\geq 0\; \text{and}\;
k\not \in [-\epsilon-\sqrt{V_--V_+}, \sqrt{V_--V_+}+\epsilon]\}$,
any $\epsilon >0$.  A similar statement
holds for $\Rp(z(-k))/\Tm(z(-k))$ for $k$
in the closed lower half-plane.
The relationships (\ref{eq:T+-}) and (\ref{eq:rt}) combine to
give, for $k\in \Real$, $k\not \in [-\sqrt{V_--V_+}, \sqrt{V_--V_+}]$,  
$$\Tm(z(-k))\Rp(z(k))+\Rm(z(-k))\Tm(z(k))=0.$$
This ensures
that $F(k)$ is holomorphic on 
$\Complex \setminus  [-\epsilon-\sqrt{V_--V_+}, \sqrt{V_--V_+}+\epsilon]$
for any $\epsilon >0$.
From Propositions
\ref{prop:t} and \ref{prop:r}, we know that
$$h_F(\varphi)= 2b_1\sin(|\varphi|).$$  Thus applying Theorem 
\ref{thm:complexanal2}, we obtain, using the notation of that theorem,
$$n_{F,\sqrt{V_- -V_+}}(r)=\frac{4b_1}{\pi}r+o(r).$$

It now remains to relate the zeros of $F$ to the poles of the resolvent
in the desired region. 
The function $\Tm$ is nonzero, except, perhaps, at points which project
to $V_+$.  
 By Proposition \ref{prop:rsm} and
(\ref{eq:r+}), the zeros of $F$ in the upper half plane correspond, with 
multiplicity, to poles of the resolvent on the sheet $\{ z\in \hat{Z}:
\Im r_+<0, \Im r_->0\}.$  Similarly, using (\ref{eq:r-}), the zeros of
$F$ in the lower half-plane correspond, with multiplicity, to poles
of the resolvent on the sheet $\{ z\in \hat{Z}:
\Im r_+>0, \Im r_-<0\}.$  By Proposition \ref{prop:rsm2}, the zeros of 
$F(k)$ with $k\in \Real$ coincide with poles at $z\in \hat{Z}$ with 
$\pro(z)\in \Real$, $\pro(z)>V_-$, and $r_+(z)$ and $r_-(z)$ having opposite
signs.
\end{proof}

Finally, we include a proposition which we shall need 
for our inverse results.
\begin{prop}
\label{prop:capp}
We have
$$\sum_{ \genfrac{}{}{0pt}{}
{z_j\in \mcr}{
 r_+(z_j)\not = 0}} \frac{|\Im (r_+(z_j))|}{|r_+(z_j)|^2} < \infty.$$
\end{prop}
\begin{proof}
 The proof follows from a application of Carleman's Theorem, using the 
variable $k=r_+(z)$, to, in turn, $\Rm$, $\Rp$, and $\Rm \Rp - \Tm\Tp$
in the physical plane.  We also use Propositions 
\ref{prop:rsm2}, \ref{prop:t}, and \ref{prop:r}.
\end{proof}

\section{Inverse Results}\label{s:inverse}
In this section we prove our inverse results.
We must recover the reflection and transmission coefficients,
$R_{\pm}$, $T_{\pm}$ on the boundary of the physical space.  Because 
$V'$ is compactly supported, this information is enough to recover the 
eigenvalues and the norming constants.  Then, using results of 
\cite{c-k}, we recover the potential.

We will work with functions meromorphic in the plane whose zeros and poles
are determined by $\mcr$.   
We can recover these functions, up to a finite number of unknown constants,
 by applying the Weierstrass
factorization theorem.  The difficulty is to find enough such functions
to be able to recover $R_{\pm}$ and $T_{\pm}$.

\begin{lemma}\label{l:inverse1}
 For a real-valued
steplike $V$, with $V'$ compactly supported,
$V_+$, $V_-$, and $\mcr$ determine
$$
 \Rm(z)\Rm(\omega_{-}(z)), \; 
\Rp(z)\Rp(\omega_+(z)), \; \text{
and }  \Tm(z)/\Tm(\omega_{+-}(z)).$$  Furthermore, $\Tm(z)\Tm(\omega_-(z))$
and $\Tm(z)\Tm(\omega_+(z))$ are determined up to constant real multiples.
\end{lemma}
\begin{proof}
Note that $\Rm(z)\Rm(\omega_-(z))$ is a meromorphic function of 
$r_+(z)$.  The poles of $\Rm(z)$ are the same as the elements of $\mcr$,
except that $\Rm(z)$ is always regular at $\pro^{-1}(V_+)$.  The zeros
of $\Rm(z)$ are those $z'\in \hat{Z}$ such that $\omega_+(z')\in \mcr$, 
$\pro(z')\not
= V_+$.  Therefore, using Proposition \ref{prop:r} and the Weierstrass  
Factorization Theorem,
$$\Rm(z)\Rm(\omega_{-}(z))= \gamma_1 e^{\delta_1 r_+(z)}
\prod_{z_j\in \mcr, \pro(z_j)\not = V_+} \frac{r_+(z_j)+r_+(z)}
{r_+(z_j)-r_+(z)}$$
where $\gamma_1$ and $\delta_1$ are constants to be determined.  The product
converges by Proposition \ref{prop:capp}.
But
$$\Rm(z)\Rm(\omega_-(z))=\frac{\Rm(z)}{\Rp(z)}(\Rm(z)\Rp(z)-\Tm(z)\Tp(z))$$
using (\ref{eq:T+-}), (\ref{eq:r+}),
(\ref{eq:rt}) and (\ref{eq:norm}).  When $z$ lies on the physical sheet,
 using Propositions \ref{prop:t} and \ref{prop:r}, this is a function of
completely regular growth of type determined by the length of the convex
hull of the support of $V'$.  By Theorem \ref{thm:asym1}, this is determined
by the resonances, so that the resonances determine $\delta_1$.  

Fix $z'\in \pro^{-1}(V_+)$.  Then $\Rm(z')$ is $1$ or $-1$, by (\ref{eq:r+}).
It is $-1$ if $z'\not \in \mcr$ and $1$ if $z'\in \mcr$.  This fixes
the value of $\Rm(z)\Rm(\omega_-(z))$ for $r_+(z)=0$, and thus determines 
$\gamma_1$.
 
A similar argument determines $\Rp(z)\Rp(\omega_+(z))$ as a meromorphic
 function
of $r_-(z)$.

Next we consider $\Tm(z)\Tm(\omega_-(z))$ and $\Tm(z)\Tm(\omega_+(z))$, which
are meromorphic functions of $r_+(z)$ and $r_-(z)$, respectively.
Note that $\Tm(z)$ is nonzero, except, perhaps, if $z\in \pro^{-1}(V_+).$
If $z'\in \pro^{-1}(V_+)$, then $\Tm$ has a zero of order $1$ at 
$z'$ if $z'\not \in \mcr$ and is nonzero at $z'$ if $z'\not \in \mcr$.
Thus 
$$\Tm(z)\Tm(\omega_-(z))= \gamma_2e^{\delta_2r_+(z)}
(r_+(z))^{\alpha_+} \prod _{z_j \in \mcr,\; \pro(z_j)\not = V_+}
\frac{1}{1-\frac{r_+(z)}{r_+(z_j)}}$$
where 
$\alpha_+=1$ if $\pro^{-1}(V_+)\cap \mcr=\emptyset$ and $\alpha_+=0$ if
$\pro^{-1}(V_+)\cap \mcr \not = \emptyset$.  (Note that (\ref{eq:r-})
ensures that at most one of the points which projects to $V_+$ can
be a resonance.)  Again, the product converges by Proposition \ref{prop:capp}.

To see that $\delta_2$ is determined by $\mcr$, note that
$$\Tm(z)\Tm(\omega_-(z))=\frac{-r_+(z)\Tm(z) \Tp(z)}{r_-(z)\Rp(z)}.$$
Taking the reciprocal, we obtain a function which is, when $z$ is in 
the physical space, a function of completely regular growth as a function
of $r_+(z)$.  Just as for 
$\Rm(z)\Rm(\omega_-(z))$, the type is determined by $\mcr$, and thus 
$\delta_2$ is determined.  Since 
$\Tm(\omega_+(z))\Tm(\omega_{+-}(z))$ is the complex conjugate of 
$\Tm(z)\Tm(\omega_-(z))$ when $V_+<\pro(z)<V_-$, $\gamma_2$ is determined
up to a real multiple.

A similar argument shows that $\Tm(z)\Tm(\omega_+(z))$ is determined by 
$\mcr$, except for a real multiple.

Now 
$$\frac{\Tm(z)}{\Tm(\omega_{+-}(z))}=\frac{\Tm(z)\Tm(\omega_+(z))}
{\Tm(\omega_+(z))\Tm(\omega_{+-}(z))}
$$
is therefore determined up to a constant, real, multiple.  But, by 
Proposition \ref{prop:t}, 
$\Tm(z)/\Tm(\omega_{+-}(z))\rightarrow 1$ when $z$ lies on the boundary
of the physical sheet and $|z|\rightarrow \infty$.  This, then, completely
determines $\Tm(z)/\Tm(\omega_{+-}(z)).$
\end{proof}

The next lemma makes an assumption, that
$\Tm(z)\Tm(\omega_-(z))$ is determined, that we have yet to prove.   Proving
this assumption is the sticking point in proving the inverse result in
general.  We are able to prove it in some special cases, and this is 
what enables us to prove the specialized inverse results.
\begin{lemma}\label{l:assumet}  If $\Tm(z)\Tm(\omega_-(z))$ is known, and 
$V\in L^{\infty}(\Real;\Real)$ is steplike, then $\mcr$, $V_+$, and 
$V_-$ determine $\Rm(z)$.
\end{lemma}
\begin{proof}
Using (\ref{eq:rt+}), 
$$\Tm(z)\Tm(\omega_-(z))=
\frac{\Tm(z)\Tm(\omega_{+-}(z))}{\Rm(\omega_{+-}(z))}.$$

For the remainder of this proof, we assume that $z$ lies on the boundary
of the physical space, with $\pro (z)>V_-$.  Then 
$$\overline{T}_-(z)=\Tm(\omega_{+-}(z)),\; \overline{R}_-(z)=
\Rm(\omega_{+-}(z)).$$
Then, away from the zeros of 
$\Rm(\omega_{+-}(z)),$ $\Tm(z)\Tm(\omega_-(z))$ determines the 
argument of $\Rm(\omega_{+-}(z))$, modulo $2\pi$.  Using (\ref{eq:norm})
and (\ref{eq:T+-}),
$$\frac{r_-(z)}{r_+(z)}\frac{\Tm(z)\Tm(\omega_{+-}(z))}{\Rm(\omega_{+-}(z))}
=
\frac{1}{\Rm(\omega_{+-}(z))}-\Rm(z).$$
Let $\rho =\rho(z)=|\Rm(z)|.$  Then
$$
f(z)=
\left|\frac{r_-(z)}{r_+(z)}\frac{\Tm(z)\Tm(\omega_{+-}(z))}{\Rm(\omega_{+-}(z))}
\right|= \frac{1}{\rho}-\rho, $$ where $f(z)\geq 0$ is known.  Thus
$$\rho =\frac{-f\pm \sqrt{f^2+4}}{2}.$$
If $\Rm(z) \not =0$ for all
 $z$ with $\pro (z) \in [V_-,\infty)$, it is clear
that we must take the ``$+$'' sign since $\Rm(z)\rightarrow 0$ as $|z|\rightarrow \infty$ (recalling that $z$ lies on the boundary of the physical sheet,
and using (\ref{eq:norm}) and Proposition \ref{prop:t}).  On the other hand, if $\Rm(z)$ does have such a zero, it is easy to see that we must take the ``$+$'' sign in the choice of $\rho$ to get $\rho =0$ at a pole of $f$.

Thus, for $z$ on the boundary of the physical space with $\pro(z)\geq V_-$, we
know the argument (modulo $2\pi$)
 and norm of $\Rm$, and thus know $\Rm$ there.  Knowing
$\Rm$ on an interval uniquely determines it on all of $\hat{Z}$.
\end{proof}

\begin{lemma}\label{l:easypart}
Suppose $V \in L^{\infty}(\Real;\Real)$ is steplike,
 $V'$ has compact support, and
$\Tm(z)\Tm(\omega_{+}(z))$ is known.  Then $\mcr$, $V_+$, and $V_-$ 
determine $\Rp(z)$ and $\Tm(z)$.
\end{lemma}
\begin{proof}
By Lemma \ref{l:inverse1}, $\Rp(z)\Rp(\omega_+(z))=\Rp(z)/\Rp(\omega_{+-}(z))$
and $\Tm(z)/\Tm(\omega_{+-}(z))$ are determined.  By Lemma \ref{l:assumet},
$\Rm(z)$ is determined, and by (\ref{eq:T+-}) and
(\ref{eq:norm}) this determines 
$\Tm(z)\Tm(\omega_{+-}(z)).$  By (\ref{eq:norm2}), $\Rp(z)\Rp(\omega_{+-}(z))$
is determined.  Therefore, $(\Tm(z))^2$ and $(\Rp(z))^2$ are fixed.  But then
$\Tm(z)$ is determined by $\Tm(z)\rightarrow 1$ as $|z|\rightarrow \infty$ with
$z$ on the boundary of the physical sheet, and $\Rp(z)$ is determined by 
$\Rp(z_1)=-1$ when $\pro (z_1)=V_-$.
\end{proof}

By Lemmas \ref{l:assumet} and \ref{l:easypart}
and applying results of \cite{c-k}, knowing $\mcr$, $V_+$, and $V_-$
will determine a steplike $V$ (with $V'$ compactly supported) {\em provided}
that $\Tm(z)\Tm(\omega_+(z))$ is fixed.  Recall by Lemma \ref{l:inverse1}
that $\Tm(z)\Tm(\omega_+(z))$ is fixed by $\mcr$
up to a real, constant, multiple. 
We show how to fix this multiple in certain cases.

\begin{lemma}  If $\pro^{-1}(V_+)\cap \mcr \not = \emptyset$, then, if 
$V\in L^{\infty}(\Real;\Real)$ is steplike with $V'$ compactly supported,
then $V_+$, $V_-$, and $\mcr$ determine $\Tm(z)\Tm(\omega_+(z))$.
\end{lemma}
\begin{proof}
By Lemma \ref{l:inverse1}, $\Tm(z)\Tm(\omega_+(z))$
and $\Tm(z)\Tm(\omega_-(z))$ are determined up to 
real constant multiples.  We fix that multiple by determining what happens
at one point.  Fix $z_0\in \pro^{-1}(V_+).$  If $z_0 \in \mcr$, then
$\Rm(z_0)=1$ and $\Rm(\omega_-(z_0))=-1$.  If $z_0 \not \in \mcr$, then
$\Rm(z_0)=-1$ and $\Rm(\omega_-(z_0))=1$.  By (\ref{eq:rdiff}), 
$$-\frac{r_-(z)\Tm(z)\Tm(\omega_-(z))}{r_+(z)}=\Rm(z)-\Rm(\omega_-(z)).$$
We thus obtain
$$\lim_{z\rightarrow z_0}-\frac{r_-(z)\Tm(z)\Tm(\omega_-(z))}{r_+(z)}
= \left\{ \begin{array}{ll}
-2 & \text{if}\; z_0 \not \in \mcr \\
2 & \text{if}\; z_0 \in \mcr
\end{array}
\right.
$$
Thus $\Tm(z)\Tm(\omega_-(z))$ is fixed, and so is $\Tm(\omega_+(z))\Tm(\omega_{+-}(z)). $
Since 
$$\frac{\Tm(z)\Tm(\omega_+(z))}{\Tm(\omega_+(z))\Tm(\omega_{+-}(z))}
$$ is fixed by Lemma \ref{l:inverse1}, $\Tm(z)\Tm(\omega_+(z))$ is determined
as well.
\end{proof}

\begin{lemma} Suppose $V$ is a real-valued steplike function, $V_+$ and $V_-$
are known, and it is known {\em a priori} that $V=p+V_{\beta}$ for some
(unknown)
$\beta\in \Real$ and (unknown) $p\in C_{\comp}^0(\Real)$
with $p'\in L^1(\Real)$.  Then
$\mcr$ determines $\Tm(z)\Tm(\omega_+(z)).$
\end{lemma}
\begin{proof}
Recall that $$\Tm(z)\Tm(\omega_+(z))=-
\frac{\Tm(z)\Tm(\omega_{+-}(z))}{\Rp(\omega_{+-}(z))}$$
and that by Lemma \ref{l:inverse1} $\Tm(z)\Tm(\omega_{+}(z))$ is determined up
to a real constant multiple.  Suppose that $z$ lies on the boundary of the
physical space.  Then, by Lemma \ref{l:Vnice}, 
$|\Rp(\omega_{+-}(z))|=\frac{|V_+-V_-|}{4(r_+(z))^2}+o((r_+(z))^{-2})$ 
as $|z|\rightarrow \infty$.  Then
$$|\Tm(z)\Tm(\omega_+(z))| =\frac{4(r_+(z))^2}{|V_+-V_-|  }+o((r_+(z))^{-2})$$
when $|z|\rightarrow \infty$ and $z$ is on the boundary of the physical sheet.
This determines the (real) constant multiple, and thus $\Tm(z)\Tm(\omega_+(z))$.
\end{proof}

\begin{lemma}
Suppose that $V$ is a real-valued steplike function, and $V_+$, $V_-$, and
$\mcr$  are 
known, with
$\pro^{-1}(V_+)\cap \mcr = \emptyset$.  
Let $z_0\in \hat{Z}$ such that $\pro (z_0)=V_+$.  Then
if $\Tp(z_0)/\Tp(\omega_-(z_0))>0$, $\Tm(z)\Tm(\omega_+(z))$ is 
determined by $\mcr$.
\end{lemma}
We remark that by Lemma \ref{l:inverse1}
and (\ref{eq:T+-}), $\Tp(z_0)/\Tp(\omega_-(z_0))$ is 
determined by $\mcr$.  This ratio must be real valued because for
$z$ with $V_+\leq \pro(z)\leq V_-$, $\Tp(z)=\overline{T}_+(\omega_+(z))$
and $\omega_+(z_0)=z_0$.
 Checking this ratio in the simple case of a 
potential which takes on only three (distinct) values shows that this
ratio can have either sign (see e.g. \cite[Appendix 1]{weder}).
\begin{proof}
Since $z_0,\; \omega_-(z_0) \not \in \mcr$, $\Rm(z_0)=-1=\Rm(\omega_-(z_0)).$
By (\ref{eq:T+-}) and (\ref{eq:rt}), and using the fact that $\omega_+(z_0)
= z_0$,
we get 
$$\Rp(z_0)=-\Rm(\omega_-(z_0))\frac{\Tp(z_0)}{\Tp(\omega_-(z_0))}.$$
Our assumption on the sign of $\Tp(z_0)/\Tp(\omega_-(z_0))$ shows that
$\Rp(z_0)>0$, and thus $\Rp(\omega_-(z_0))>0$ as well.

Let $z_1 \in \hat{Z}$ have $\pro (z_1)=V_-$.  We have already noted that 
since $z_1\not \in \mcr$, $\Rp(z_1)=-1$.  
Because $\Rp$ is continuous at all points $z$ with $\pro(z)\in [V_+, V_-]$
and is nonzero is this range, there is some point $z'$ with $\pro(z')
\in (V_+,V_-)$ and $\Rp(z')$ pure imaginary.  At this point, $\Rp(w_+(z'))=
-\Rp(z')$ and 
$r_+(z')\Tp(z')\Tp(\omega_+(z'))=-2r_-(z')\Rp(z')$ by (\ref{eq:rdiff2}).
That is, $|r_+(z)\Tp(z)T_+(\omega_+(z)/r_-(z)\Rp(z)|$ is a maximum at $z=z'$ for 
$\pro(z)\in [V_+,V_-]$  Note that for $\pro(z)\in 
[V_+,V_-]$, $|\Rp(z)|=(\Rp(z)\Rp(\omega_+(z)))^{1/2}$ is known.  Thus the 
norm of $r_+(z')\Tp(z')\Tp(\omega_+(z'))/r_-(z')=-2\Rp(z')$ is known, and since
$\Tp(z)\Tp(\omega_+(z))$ was already determined up to a real constant 
multiple, it is completely determined.
\end{proof}

We summarize our inverse results in the following theorem.
\begin{thm}  Knowledge of $V_+$, $V_-$, and $\mcr$
uniquely  determines, up to translation, a real
steplike
potential $V$ with $V'$ compactly supported provided at least one of the
following criteria is met:
\begin{itemize}
\item $\Pi^{-1}(V_+)\not = \emptyset$
\item It is known {\em a priori} that $V=V_{\beta}+p$ for some (unknown)
$\beta\in \Real$ and (unknown) $p\in C^0_{\comp}(\Real)$ with
$p'\in L^1(\Real)$.
\item For $z_0\in \Pi^{-1}(V_+),$ $\Tp(z_0)/\Tp(\omega_-(z_0))>0$.
\end{itemize} 
\end{thm}

\small
\noindent
{\sc 
Department of Mathematics,
University of Missouri,
Columbia, Missouri 65211\\
{\tt tjc@math.missouri.edu}

\begin{thebibliography}{99}

\bibitem{b-k-w} B.M. Brown, I. Knowles, and R. Weikard, 
{\em On the inverse resonance problem}, preprint.

\bibitem{ace}
T. Christiansen, {\em Scattering theory for manifolds with asymptotically 
cylindrical ends}. Journal of Functional Analysis {\bf 131}, 2, 
(1995), 499-530.

\bibitem{cep} T. Christiansen,
{\em Some upper bounds on the number of resonances for manifolds
with infinite cylindrical ends}, Annales Henri Poincar\'e {\bf 3} No. 5
(2002), 895-920.


\bibitem{c-k} A. Cohen and T. Kappeler, {\em Scattering and inverse scattering
for steplike potentials in the Schr\"odinger equation}, 
 Indiana Univ. Math. J. {\bf 34}
(1985), no. 1, 127--180.


\bibitem{froese} R. Froese, {\em Asymptotic distribution of resonances
in one dimension}, J.
Differential Equations {\bf 137} (1997), no. 2, 251--272.

\bibitem{froese2} R. Froese, {\em Upper bounds for the resonance counting 
function of Schr\"odinger operators in odd dimensions},
 Canad. J. Math. {\bf 50} (1998),
no. 3, 538--546.

\bibitem{guce} L. Guillop\'{e}, {\em Th\'{e}orie spectrale de quelques 
vari\'{e}t\'{e}s \`{a} bouts}, Ann. Scient. Ec. Norm. Sup. {\bf 22},
4, (1989), 137-160.



\bibitem{k-ir} E. Korotyaev, {\em Inverse resonance scattering on the 
real line}, preprint.

\bibitem{k-hl} E. Korotyaev, {\em Inverse resonance scattering on the 
half line}, preprint.

\bibitem{levin} B. Ja. Levin, {\em Distribution of zeros of entire functions},
American Mathematical Society, Providence, R.I. 1964 viii+493 pp. 


\bibitem{tapsit} R.B. Melrose, {\em The Atiyah-Patodi-Singer Index 
Theorem}, A.K. Peters, Wellesley, MA 1993.


\bibitem{pac} L. Parnovski, {\em Spectral asymptotics of the Laplace 
operator on surfaces with cusps}, Math. Ann. {\bf 303} (1995), 281-296.


\bibitem{regge} T. Regge, {\em Analytic properties of the scattering 
matrix}, Nuovo Cimento {\bf 8} (5), (1958), 671-679.

\bibitem{sjrfb} J. Sj\"ostrand, {\em Resonances for bottles and trace 
formulae}, Math. Nachr. {\bf 221} (2001), 95--149.




\bibitem{vdegen} G. Vodev, {\em Asymptotics on the number of scattering poles 
for degenerate perturbations of the Laplacian}, J. Funct. Anal. 
{\bf 138} (1996), 295-310.

\bibitem{vodevsurvey}
G. Vodev, {\em Resonances in the Euclidean scattering}, Cubo
Matem\'atica Educacional {\bf 3} no. 1 (2001), 317-360.

\bibitem{weder} R. Weder, {\em Spectral and scattering theory
for wave propagation in perturbed stratified media}, Springer-Verlag,
New York, 1991.


\bibitem{zw1} M. Zworski, {\em Distribution of poles for
scattering on the real line}, J. Funct. Anal. {\bf 73} (2) (1987),
277-296.

\bibitem{zwrp} M. Zworski, {\em Sharp polynomial bounds on the number of 
scattering poles of radial potentials}, J. Funct. Anal. {\bf 82} (1989), 
370-403.


\bibitem{zwsurvey} M. Zworski, {\em Counting scattering poles}.  In:
 Spectral and scattering theory (Sanda, 1992), 301--331,
 Lecture Notes in Pure and Appl. Math., 161,
Dekker, New York, 1994. 


\bibitem{ZwS} M. Zworski, {\em Resonances in physics and geometry.}
Notices Amer. Math. Soc. {\bf 46} (1999), 319--328.

\bibitem{zwinverse} M. Zworski, {\em 
A remark on isopolar potentials}, SIAM J. Math. Analysis, 
{\bf 82} (6) (2001), 1823-1826. 

\end{thebibliography}
\end{document}